\documentclass[12pt]{article}
\usepackage{amssymb,amsmath}

\newtheorem{thm}{Theorem}

\newtheorem{lem}[thm]{Lemma}

\newtheorem{rem}[thm]{Remark}

\def\Z{\mathbb{Z}}

\title{Factorization and malleability of RSA modules, and counting points on elliptic curves modulo $N$ }

\author{Luis V. Dieulefait, Jorge  Urroz\\ 
 Dept. mat. i inform\`atica, Universitat de Barcelona, \\
 Dept math. UPC
 }

\date{}

\begin{document}\maketitle

\begin{abstract} In this paper we address two different problems related with the factorization of an RSA module $N$. First we can show that factoring is equivalent in deterministic polynomial time to counting points on a pair of twisted Elliptic curves modulo $N$. Also we settle the malleability of factoring an RSA module, as described in \cite{pavi}, using the number of points of  a single elliptic curve modulo $N$, and Coppersmith's algorithm. 
\end{abstract}
\section{Introduction}

The motivation of this note is twofolded. First we address the problem of malleability  of an RSA module $N$ and, from there, the equivalence between factoring and counting the number of points on an elliptic curve modulo $N$.

\

Malleability started in the paper by Pailler and Villar \cite{pavi}, while studying the existence of a tradeoff between one-wayness and chosen ciphertext security,  already observed  back to the eighties for example  in \cite{rabin,williams,gomiri}. In some sense, one cannot achieve one-way encryption with a level of security equivalent to solve certain difficult problem, at the same time as the cryptosystem being IND-CCA secure with respect to it.

\

Even though this paradox has been observed, it has not been formally proved except in the case of factoring-based cryptosystems in which Pailler and Villar \cite{pavi} clarified the question reformulating the paradox
 in terms of key preserving black-box reductions and proved that if factoring can be reduced in the standard model to breaking one-wayness of the cryptosystem then it is impossible to achieve chosen-cyphertext security.

\

After this they introduce the notion of malleability of a key generator  and, with it,  they are able to extend the result from key preserving black box reductions to the case of arbitrary black box reductions.
This notion, which we give explicitly in Section \ref{sec:mal}, captures a very basic fact in arithmetic: intuitively, one tends to believe that the problem of factoring a given number $n$ (an RSA modulus) is not made easier if we know how to factor other numbers $n'$ relatively prime to $n$. If this is true we say that  factoring is non-malleable.

\

As the authors themselves stress in \cite{pavi}, it is very important to study non-malleability of key generators and, in fact, they conjecture that most instance generators are non-malleable, although no arguments are given to support this belief.

\

In \cite{diu1} we address this question and notice that the freedom of selecting the new number $n'$ breaks the independent behaviour of prime numbers, and hence we produce an explicit $n'$ which makes factorization malleable. In other words, given any RSA modulus $n$  we  prove the existence of a polynomial time reduction algorithm from factoring $n$ to factoring certain explicit numbers $n'$, all relatively prime to $n$.

\

The numbers given in \cite{diu1} are very simple: given the RSA modulus $n=pq$, the factorization of $n'=m^n-1$, where $m$ is a primitive root of the smallest prime dividing $n$, allows us  to factor $n$ in polynomial time. However, this does not give a complete satisfactory answer for several reasons. First one could think of $n'$ to be of exponential size and then out of the scope of the question. However, as we mention in \cite{diu1}, one can think of $n'$ as a collection of exactly $n$ ones when it is written in $m$-ary, and we just need the factors of $n'$ modulo $n$, a data that has the same size as the given number. In any case, it still persists the restlessness of knowing whether or not in an small interval centered in $n$ we can find an explicit $n'$ which can help to factor $n$.

\


\

In this paper we address this question, given an affirmative answer to the  malleability of the problem of factoring by showing a number of the same size of $N$ whose factorization allows us to factor $N$ with a  algorithm that runs in polynomial time.

\

The notion of malleability  rests in measuring the difference between suitable  Game $0$ and Game $1$ as defined in \cite{pavi}, more precisely: if we compare the success probability of factoring $N$ with an Oracle which can solve any problem that can be reduced to factoring (Game $0$) with the one when using an Oracle having the extra ability of factoring numbers which are relatively prime to $N$ (Game $1$) the probability of factoring $N$ increases significantly (see \cite{pavi}, Section 4.1,  for more details).

\

Concretely we will prove that given a random elliptic curve $E$ defined modulo $N$, where $N$ is an RSA module, and assuming that its number of points $E_N$ is known, by further knowing the factorization of $E_N$ we can produce a deterministic polynomial time algorithm that factors $N$. The key tool in our proof will be the result of Coppersmith (see \cite{coperbi}) that allows to factor an integer by knowing only certain bits of one of its prime factors. It is worth remarking that it is not known how to factor $N$ only with the number $E_N$ as input.

\

While proving this statement  another interesting problem treated widely in the literature, (see \cite{vimoma} and \cite{kuko} for related results) showed up:

\

{\bf Problem} Is  factoring  $N$ equivalent to counting the number of points of elliptic curves modulo $N$.?

\

We were  lucky and we could also give a definite answer  by proving the following theorem:

\begin{thm}\label{teo:factor} Given the number of points, affine or projective, of any elliptic curve and one of its twists modulo $N$   we can factor $N$ in deterministic polynomial time.
\end{thm}

\

\begin{rem} As we said this  problem has been addressed in  \cite{kuko}. We should stress that their results are based in an assumption on the distribution of the number of points on elliptic curves over finite fields which is not accurate. But more than that,  the reduction algorithm from counting the number of the elliptic curve modulo $N$ to factoring $N$ in their case is probabilistic while here it is proved to be deterministic. Moreover, in terms of malleability, what we do in Section \ref{sec:mal} involves taking a single elliptic curve, and suceeds with probability $1$, while the results in \cite{kuko} need to consider many elliptic curves to have positive probability to factor $N$. Finally  the method used in that paper only works for the number of projective points on the elliptic curve, not covering the affine case as we do.

\end{rem}

\

The structure of the paper goes as follows: In Section   \ref{sec:fact} we prove Theorem \ref{teo:factor}, while Section \ref{sec:mal} is dedicated to the problem of malleability of factoring.

\

\section{Factorization}\label{sec:fact} 

\

Let $N\in \Z$. Given an elliptic curve $E:=\{y^2=x^3+ax+b\}$ over $\Z/N\Z$,   we will denote by $E_d$ its quadratic twist $E:=\{dy^2=x^3+ax+b\}$. $E(\Z/N\Z)$ will be the number of points of $E\pmod N$, and $E(\Z/N\Z)^*$ the number of affine points of the curve $\pmod N$. We know that if $N=pq$, then  $E(\Z/N\Z)=E(\Z/p\Z)\times E(\Z/q\Z)$, and that for any prime $l$,  $E(\Z/l\Z)=l+1-a_l$, where $a_l$ is the trace of the Frobenius endomorphism of $E\pmod l$.

\

\noindent  Let $N=pq$ be an RSA modulus, and $d$ an integer such that $\left(\frac dp\right)=-1$ or $\left(\frac dq\right)=-1$. We will use the abuse of notation  $E=E(\Z/N\Z)$ (or $E=E(\Z/N\Z)^*$), given by
$$
E=(P-a_p)(Q-a_q)=PQ-Pa_q-Qa_p+a_pa_q,
$$
where $P=p+1, Q=q+1$ in the projective case and $P=p, Q=q$ in the affine case.

\

There are three options for $E_d(\Z/N\Z)$ (or $E_d(\Z/N\Z)^*$), which will be denoted by $\hat E,\tilde E,\bar E$ respectively, depending on the Legendre simbols  $\left(\frac dp\right)$ and $\left(\frac dq\right)$,
\begin{eqnarray*}
&&\hat E=(P+a_p)(Q+a_q)=PQ + Pa_q+Qa_p+a_pa_q,\\
&&\tilde E=(P-a_p)(Q+a_q)=PQ +Pa_q-Qa_p-a_pa_q, \\
&&\bar E=(P+a_p)(Q-a_q)=PQ -Pa_q+Qa_p-a_pa_q.
\end{eqnarray*}

Then,
\begin{equation}\label{eq:sum}
E+\hat E+\tilde E+\bar E=4PQ,
\end{equation}
while
$$
4PQ=\frac{(E+\tilde E)(E+\bar E)}E=E+\tilde E+\bar E+\frac{(\tilde E\bar E)}E=4PQ-\hat E+\frac{(\tilde E\bar E)}E,
$$
so
\begin{equation}\label{eq:product}
E\hat E=\tilde E\bar E.
\end{equation}

\begin{lem}\label{lem:twofour} Knowing two among $\tilde E,\hat E,\bar E$ and $E$, we know the four of them.
\end{lem}\
{\bf Proof.}  We split the proof in $2$ cases.

\medskip

\noindent {\bf Case 1.}
We suppose $E$ and $\hat E$ are known. The case in which $\tilde E$ and $\bar E$ are known is analogous. Then we compute its product, $M=E\hat E$ and its sum $L=E+\hat E$, and we have
\begin{eqnarray*}
&&\tilde E\bar E=M\\
&&\tilde E+\bar E=4PQ-L
\end{eqnarray*}
so $\tilde E$ and $\bar E$ are the solutions of the quadratic polynomial $X^2-(4PQ-L)X+M.$
$$
\tilde E=\frac{4PQ-L+\sqrt{(4PQ-L)^2-4M}}{2},\qquad \bar E=\frac{4PQ-L-\sqrt{(4PQ-L)^2-4M}}{2}
$$
\medskip

\noindent {\bf Case 2.} Suppose $E$ and $\tilde E$ are known. The cases in which the pairs $(E,\bar E)$, $(\hat E,\tilde E)$ and $(\hat E,\bar E)$ are known, are analogous. Then, compute the quotient $\frac{E}{\tilde E}=M$ and the sum $E+\tilde E=L$. Hence, $\frac {\bar E}{\hat E}=M$, or $\bar E=M\hat E$, and  by (\ref{eq:sum})
$(M+1)\hat E=4PQ-L,$ or
$$
\hat E=\frac{4PQ-L}{(M+1)},\qquad \bar E=\frac{M(4PQ-L)}{(M+1)}.
$$

\begin{thm}\label{teo:twist} Knowing either $E(\Z/N\Z)$ and $E_d(\Z/N\Z)$ or $E(\Z/N\Z)^*$ and $E_d(\Z/N\Z)^*$
for $\left(\frac dp\right)=-1$,  we can factor $N$ in polynomial time.
\end{thm}
In the projective case, we compute the four integers $E, \hat E,\tilde E,\bar E$ by Lemma \ref{lem:twofour} and then its sum to compute $PQ$. With $PQ$ and $N$ we factor $N$.

\

In the affine case, we again compute the four integers $E, \hat E,\tilde E,\bar E$, and then note that
 $E+\tilde E=2q(p-a_p)$, has $q$ as a common factor with $N$, so computing the gcd with $N$, we factor $N$.

 \

In both cases observe that, in principle, we do not know which one is $E_d(\Z/N\Z)$, so we will have to make two computations.

\begin{thm} Under ERH factoring an RSA modulus $N=pq$ is polynomial time equivalent to count the number of points, affine or projective, of any elliptic curve $E$ modulo $N$.
\end{thm}

\medskip

{\bf Proof.} Let $E$ be an elliptic curve. Then, knowing the factorization of $N$, we can compute $E_N$ by Schoof's algorithm \cite{schoof}.

\

Now suppose we know $E_N$. from  \cite{ankeny},  we know that, under ERH, the smallest quadratic nonresidue modulo $p$, call it $d$, is of size $O((\log p)^2)$. Hence, apply the previous Theorem  \ref{teo:twist} to the pair $E,E_d$.

\

Recall that, as of today, we can compute the number of points modulo $N$ by baby step giant step, since $E\pmod N$ has group structure, in $O(N^{1/4+\varepsilon})$ which is exponential.

\
\section{Malleability}\label{sec:mal}

As in previous sections, let $N=pq$ be an RSA module.  We recall that in order to prove that factoring is malleable we need to find a number relatively prime to $N$ and of the same size,  which factorization will allow us to factor $N$ in deterministic polynomial time. We consider a random elliptic curve $E\pmod N$, and we let $E_N$  be the number of affine points, while $E^*_N$ will be the number of points including the points at infinity. We can assume that we have at our disposal an Oracle that computes any of these two numbers. Since this computation can be reduced to the factorization of $N$ thanks to Schoof's algorithm, this corresponds to Game 0 in the setup described in the introduction while defining malleability. Note also that, as we have already observed, there is no known polynomial time algorithm that can factor $N$ using this information.

\

Assume now that we have access to an auxiliary Oracle that can factor any number relatively prime to $N$. Using it, we factor the number $E_N$ (or $E^*_N$) and we will show in the following theorem that from this we can factor $N$ in polynomial time, thus concluding that Game $1$ has solved the factorization problem that was not achieved by Game 0, which shows that factorization of RSA modules is malleable.

\

\begin{thm} Given $N=pq$  where $p,q$ are prime numbers, and an elliptic curve $E\pmod N$,  there exists a polynomial time algorithm in $\log N$ such that  with input $N$ and the factorization of  $E_N$ or $E^*_N$, it outputs the  factors of $N$, $p$ and $q$ with probability one.
\end{thm}

{\bf Proof.} As we mentioned in the introduction, we will use a well known result of Copersmith, which allows us to find a factor of an integer by just knowing certain part of its highest bits. For convenience we include this result now
\begin{thm} (Coppersmith) If we know an integer $N=pq$ and we know the high order
$(1/4)(log_2N$  bits of $p$, then in timee polynomial in $log N$ we ran discover $p$
and $q$.
\end{thm}

Observe that it would be sufficient by knowing $(1/4)(log_2N-O(\log\log N$, since we could try the unknown bits up to $(1/4)(log_2N$ one by one in polynomial time.

\

Now,  recall again that, by Hasse's theorem the factor found $q-a_q+1$ is at distance $|a_q-1|\le 2\sqrt q+1\le 2N^{1/4}+1$ of $q$ which is a factor of $N$. By bounding the distance, we know certain of the highest bits of $q$ from those of $q-a_q+1$. 
 
 \
 
 In particular, let us suppose that two integers $x<y$ are at distance $y-x=2^t+R$ where $R<2^t$. We can write $x=M_x2^t+R_x$, $y=M_y2^t+R_y$ with $R_x<2^t$, $R_y<2^t$ and $-2^t< R_y-R_x<2^t$. Then, $y-x=(M_y-M_x)2^t+R_y-R_x$, which gives $M_y=M_x+1$ or $M_y=M_x+2$ and, hence, from the highest bits of $y$ up to $t$ of $x$ we know those of $y$ and viceversa. 
 
 \
 
In our case, the distance is bounded by $2q^{1/2}$ so we know up to $t=[\log_2 q/2]+1$ of the highest bits of $q$. By division, we also know  up to $t$ of the highest bits of $p$. But $\sqrt N\le p=M_p2^t+R_p\le (M_p+1)2^t$, and so $M_p\ge \sqrt N/(2^t+1)\ge N^{1/4}$ and, hence, we can apply Coppersmith algorithm to find $p$, thus factoring $N$. 

\

Let us stress that having factored the number $E^*_N$, from the average of the divisor function
$$
\sum_{n\le x}d(n)=x\log x+(2\gamma-1)x+o(x)
$$
we deduce that, with probability $1$,  the number of factors of an integer is of the order of the logarithm of it. Hence, once the auxiliary Oracle gives the factors of $E_N$, we will apply Coppersmith one by one finding $q$  in polynomial time in log $N$.

\begin{rem}The case in which $E_N$ is given, is similar  and we leave the details to the reader.
\end{rem}

\

\subsection{Small difference}

Even though malleability is fully proved in the last section, we include this section as a small remark in the negligibale  case in which $E_N$ and $E_N^*$ have an exponential number of divisors, but the two prime factors $p,q$ are not too far from each other.

\

In order to construct a module RSA we tipically search for a couple of prime factors of the same number of bits, i.e. $q<p<2q$. However, if the two primes are very close to each other, the scheme is easy to break since the module can be factored in polynomial time. Indeed,  It is well known that if $\Delta=|p-q|<N^{1/4}$  Fermat's factorization algorithm enables to find both factors of $N$ in polynomial time and there has been an effort of the community  to improve the  exponent $1/4$ in $\Delta$  for the factorization of $N$. It is worth to mention that if the objective is  breaking the RSA scheme, rather than factoring the modulus, then the exponent can be increased all the way up to basically $1$ by means of an improved version of the attacks done by Wiener or Boneh and Durfee, (see \cite{weger}). However, for the factorization of $N$ not too much more is known.
In \cite{erragre} the authors claim in an apparently unpublished work that  are able to factor an RSA module $N=pq$ even when the difference is of order $|p-q|<N^{1/3}$. 

\

We devote this section to recover  $\Delta<N^{1/3}$ using malleability techniques: in particular the factorization of the number of points of a random elliptic curve modulo $N$, together with  a simple application of an argument of elementary geometry  attributed to Heron of Alexandria which says that  in any triangle, the product of the length of its three sides equals four times the area times the radius of the circumscribed circle. 
We will assume from now that $\Delta=|p-q|<c'N^{1/3}$ for some suitable constant $c'$.

\

In our case given three points  $(x_0,y_0),(x_1,y_1),(x_2,y_2)$ of integer coordinates in the hyperbola $xy=E_N^*$, we see that the radius of the circumscribed circle is 
$$
R=\frac{((x_0x_1)^2+(E_N^*)^2)((x_0x_2)^2+(E_N^*)^2)((x_2x_1)^2+(E_N^*)^2)}{4(E_N^*)^2(x_0x_1x_2)^2},
$$
and taking $\max\{x_0,x_1,x_2\}\le \sqrt {E_N^*}$, we get 
$$
R\ge \frac{E_N^*}{4}.
$$
On the other hand, by Hasse's theorem
$$
E_N^*\ge (\sqrt p-1)^2(\sqrt q-1)^2
$$
and 
$$
(\sqrt p-1)(\sqrt q-1)=\sqrt N-\sqrt p-\sqrt q+1\ge \sqrt N/4,
$$
for $N$ sufficiently large, and so 
$$
R\ge \frac{N}{64}.
$$
Hence, by Heron of Alexandria's theorem, in an arc of the hyperbola $xy=E_N^*$ of lenght less than $(N/32)^{1/3}$ we can only have two points of integer coordinates.


\


Now recall that the lenght $L$ of an arc  of the hyperbola $xy=T$ with $a\le x\le b$ is given by 
\begin{eqnarray*}
L&=&\int_a^b\sqrt{1+\frac{T^2}{t^4}}dt\le \int_a^b1+\frac{T^2}{t^4}dt=(b-a)+\frac{T^2}3\left(\frac1{a^3}-\frac1{b^3}\right)=\\
&=&(b-a)+\frac{(b-a)T^2}3\left(\frac{b^2+ab+a^2}{(ab)^3}\right)\le (b-a)+\frac{(b-a)T^2}{a^3b}.
\end{eqnarray*}
Consider $T=E_N^*$, $b\ge N^{1/2}$ and $b-a\le cN^{1/3}$ for suitable $c$. 
Hence, by noting that 
$$
E_N^*\le (\sqrt p+1)^2(\sqrt q+1)^2=(\sqrt N+\sqrt p+\sqrt q+1)\le N+7 N^{3/4},
$$
since the primes are very close, we get from a simple computation  that the arc on the hyperbola has lenght 
$$
L\le 3cN^{1/3},
$$
 In particular, we can select $c$ small enought so $L\le (N/32)^{1/3}$ and hence it can only have at most two points of integral coordinates. Hence we can ask  the auxiliary Oracle to factor $E^*_N$ and output the at most two  factors of it lying in the interval $[a,b]$. This Oracle will give back the factor $a\le N^{1/2}-c'N^{1/3}\le q-a_q+1\le N^{1/2}\le b$, for $a_q\ge 0$ or 
$a\le N^{1/2}\le q-a_q+1\le N^{1/2}+c'N^{1/3}\le b$ if $a_q\le 0$ for some $c'\le c$. In practice $c=\frac1{3(32)^{1/3}}$ and $c'=\frac1{6(32)^{1/3}}$ are enough. Using it we can factor $N$ with Coppersmith's algorithm, as we did in the previous subsection.

 \ 
 
\begin{rem} Again the case in which $E_N$ is given, is similar and we leave the details to the reader.
\end{rem}


 \bibliographystyle{plain}
\bibliography{modularfactor}

\begin{thebibliography}{10}

\bibitem{ankeny}
N.~C. Ankeny.
\newblock The least quadratic non residue.
\newblock {\em Ann. of Math. (2)}, 55:65--72, 1952.

\bibitem{coperbi}
Don Coppersmith.
\newblock Finding a small root of a bivariate integer equation; factoring with
  high bits known.
\newblock {\em Advances in cryptology, EUROCRYPT '96, Lecture Notes in Comput.
  Sci.}, 1070:178--189, 1996.

\bibitem{weger}
B.~de~Weger.
\newblock Cryptanalysis of rsa with small prime difference.
\newblock {\em Appl. Algebra Engrg. Comm. Comput.}, 13(1):17--28, 2002.

\bibitem{diu1}
L.~Dieulefait and J.~Jim\'{e}nez~Urroz.
\newblock Small primitive roots and malleability of {RSA} moduli.
\newblock {\em J. Comb. Number Theory}, 2(2):171--179, 2010.

\bibitem{erragre}
R~Erra and C~Grenier.
\newblock The fermat factorization method revisited.
\newblock {\em https://eprint.iacr.org/2009/318.pdf}, 2009.

\bibitem{gomiri}
S.~Goldwasser, S.~Micali, and R.~L. Rivest.
\newblock A digital signature scheme secure against adaptive chosen-message
  attacks.
\newblock {\em SIAM J. Comput.}, 17(2):281--308, 1988.
\newblock Special issue on cryptography.

\bibitem{kuko}
N.~Kunihiro and K.~Koyama.
\newblock Equivalence of counting the number of points on elliptic curve over
  the ring $z_n$ and factoring $n$.
\newblock {\em Lecture Notes in Comput. Sci.}, 1043:47--58, 1998.

\bibitem{vimoma}
S.~Mart\'{\i}n, P.~Morillo, and J.~L. Villar.
\newblock Computing the order of points on an elliptic curve modulo n is as
  difficult as factoring n.
\newblock {\em Appl. Math. Lett.}, 14(3):341--346, 2001.

\bibitem{pavi}
P.~Paillier and J.~L. Villar.
\newblock Trading one-wayness against chosen-ciphertext security in
  factoring-based encryption.
\newblock In {\em Advances in cryptology---{ASIACRYPT} 2006}, volume 4284 of
  {\em Lecture Notes in Comput. Sci.}, pages 252--266. Springer, Berlin, 2006.

\bibitem{rabin}
M.~O. Rabin.
\newblock Digital signatures and public key functions as intractable as
  factorization.
\newblock {\em Technical Report MIT/LCS/TR-212}, 1979.

\bibitem{schoof}
R.~Schoof.
\newblock Elliptic curves over finite fields and the computation of square
  roots$\operatorname{mod} p$.
\newblock {\em Mathematics of Computation}, 44(170):483--494, 1985.

\bibitem{williams}
H.~C. Williams.
\newblock A modification of the {RSA} public-key encryption procedure.
\newblock {\em IEEE Trans. Inform. Theory}, 26(6):726--729, 1980.

\end{thebibliography}

\end{document}